\documentclass[12pt,reqno]{amsart}
\usepackage{fullpage, times, amssymb, amsthm, amsmath, amsfonts}
\usepackage{graphicx}
\usepackage{times}
\usepackage{enumitem}
\usepackage{mathrsfs}
\usepackage[all,cmtip]{xy}
\usepackage[usenames,dvipsnames]{xcolor}
\usepackage[margin=1.5cm]{geometry}
\usepackage{hyperref}
\hypersetup{colorlinks=true, citecolor=blue, linkcolor=blue, urlcolor=blue, pdfstartview=FitH, pdfauthor=Fromm-Goldmakher, pdftitle=Improving Burgess}

\usepackage{tikz}
\usetikzlibrary{matrix,arrows,decorations.pathmorphing}

\begin{document}

\pagestyle{empty}

\parskip0pt
\parindent10pt

\newenvironment{answer}{\color{Blue}}{\color{Black}}
\newenvironment{exercise}{\color{Blue}\begin{exr}}{\end{exr}\color{Black}}

\newtheorem{theorem}{Theorem}[section]
\newtheorem{prop}[theorem]{Proposition}
\newtheorem{lemma}[theorem]{Lemma}
\newtheorem{cor}[theorem]{Corollary}
\newtheorem{conj}[theorem]{Conjecture}

\newtheorem{exr}{Exercise}
\newtheorem{example}{Example}

\newtheorem*{defn}{Definition}
\newtheorem*{rmk}{Remark}
\newtheorem*{question}{Question}
\newtheorem*{claim}{Claim}

\newtheorem*{BigThm}{Theorem A}
\newtheorem*{MainLem}{Lemma B}

\makeatletter
\newcommand{\oset}[3][0ex]{%
  \mathrel{\mathop{#3}\limits^{
    \vbox to#1{\kern-2\ex@
    \hbox{$\scriptstyle#2$}\vss}}}}
\makeatother

\renewcommand{\mod}[1]{{\ifmmode\text{\rm\ (mod~$#1$)}\else\discretionary{}{}{\hbox{ }}\rm(mod~$#1$)\fi}}

\newcommand{\act}[2]{{\ifmmode\text{$#1$\ \rotatebox[origin=c]{-90}{$\circlearrowright$}\ $#2$}}}

\newcommand{\leg}[2]{\left(\frac{#1}{#2}\right)}

\newcommand{\dtr}{\oset[-.7ex]{\cdot}{-}}
\newcommand{\isom}{\oset
    {\widetilde{\phantom{\longrightarrow}}}
    {\longrightarrow}}

\newcommand*{\dblast}
    {\raisebox{-0.1ex}{*}\llap{\raisebox{-1.25ex}{*}}}

\newcommand{\ns}{\mathrel{\unlhd}}
\newcommand{\wh}[1]{\widehat{#1}}
\newcommand{\wt}[1]{\widetilde{#1}}
\newcommand{\floor}[1]{\left\lfloor#1\right\rfloor}
\newcommand{\abs}[1]{\left|#1\right|}
\newcommand{\ds}{\displaystyle}
\newcommand{\nn}{\nonumber}
\newcommand{\im}{\textup{im }}
\renewcommand{\ker}{\textup{ker }}
\renewcommand{\Im}{\textup{Im }}
\renewcommand{\Re}{\textup{Re }}
\renewcommand{\ni}{\noindent}
\renewcommand{\bar}{\overline}

\newcommand{\mattwo}[4]{
\begin{pmatrix} #1 & #2 \\ #3 & #4 \end{pmatrix}
}

\newcommand{\vtwo}[2]{
\begin{pmatrix} #1 \\ #2 \end{pmatrix}
}

\newsymbol\dnd 232D

\newcommand{\one}{{\rm 1\hspace*{-0.4ex} \rule{0.1ex}{1.52ex}\hspace*{0.2ex}}}

\renewcommand{\v}{\vec{v}}
\newcommand{\w}{\vec{w}}

\newcommand{\Z}{\mathbb Z}
\newcommand{\Q}{\mathbb Q}
\newcommand{\N}{\mathbb N}
\newcommand{\R}{\mathbb R}
\newcommand{\C}{\mathbb C}
\renewcommand{\L}{\mathcal L}
\newcommand{\M}{\mathcal M}
\renewcommand{\P}{\mathcal P}
\renewcommand{\a}{\alpha}
\renewcommand{\b}{\beta}
\renewcommand{\d}{\delta}
\newcommand{\e}{\epsilon}
\newcommand{\F}{\mathcal F}
\newcommand{\G}{\mathcal G}
\newcommand{\E}{\mathcal E}
\newcommand{\V}{\mathcal V}

\newcommand{\ignore}[1]{}

\newcommand*\circled[1]{\tikz[baseline=(char.base)]{
            \node[shape=circle,draw,inner sep=2pt] (char) {#1};}}

\newcommand*\squared[1]{\tikz[baseline=(char.base)]{
            \node[shape=rectangle,draw,inner sep=2pt] (char) {#1};}}

\title{Improving the Burgess bound via P\'olya-Vinogradov}
\author{Elijah Fromm}
\address{Dept of Mathematics and Statistics \\
Williams College \\
Williamstown, MA, USA  01267}
\email{emf1@williams.edu}

\author{Leo Goldmakher}
\address{Dept of Mathematics and Statistics \\
Williams College \\
Williamstown, MA, USA  01267}
\email{Leo.Goldmakher@williams.edu}



\thanks{LG is partially funded by an NSA Young Investigator grant.}


\begin{abstract}
We show that even mild improvements of the P\'olya-Vinogradov inequality would imply significant improvements of Burgess' bound on character sums. Our main ingredients are a lower bound on certain types of character sums (coming from works of the second author joint with J. Bober and Y. Lamzouri) and a quantitative relationship between the mean and the logarithmic mean of a completely multiplicative function.
\end{abstract}

\maketitle
\thispagestyle{empty}

\numberwithin{equation}{section}

\section{Introduction}

Let
$\ds
S_\chi(t) := \sum_{n \leq t} \chi(n) ,
$
where $\chi \mod q$ is a Dirichlet character.
There are two famous upper bounds on this quantity.
The first, discovered independently by P\'olya and Vinogradov a century ago, asserts
\begin{equation}
\tag{$*$} \label{eq:PV}
|S_\chi(t)| \ll \sqrt{q} \log q
\end{equation}
for any primitive $\chi \mod q$.
In particular, this implies that $S_\chi(t) = o(t)$ for all $t > q^{1/2+\epsilon}$.
Sixty years ago, Burgess \cite{Burgess} found a way to increase the range of $t$ in which $S_\chi(t)$ is small. Combining his work with a clever observation of Hildebrand \cite{Hild}, it can be shown that for all primitive real quadratic characters $\xi \mod p$,
\begin{equation}
\tag{$\dag$} \label{eq:Burgess}
\phantom{\qquad \forall t > p^{1/4-o(1)}}
S_\xi(t) = o(t)
\qquad \forall t > p^{1/4-o(1)} .
\end{equation}
Here $o(1)$ is a positive quantity which tends to $0$ as $p \to \infty$. (Burgess' bound holds for more general characters as well, but for simplicity we restrict ourselves to the special case of quadratic characters of prime conductor.)

Although neither \eqref{eq:PV} nor \eqref{eq:Burgess} has been improved in general, probably neither one is optimal.
It is believed that P\'olya-Vinogradov can be improved to
\[
\tag{{\dblast}}
|S_\chi(t)| \leq
    \big(C + o(1)\big) \sqrt q \log \log q
\]
with $C = \frac{e^\gamma}{\pi}$ for odd $\chi$ and
$C = \frac{e^\gamma}{\pi\sqrt 3}$ for even $\chi$. Indeed, Montgomery and Vaughan \cite{MV} proved that GRH implies (\dblast) for some constant $C$; more recently, Granville and Soundararajan \cite{GSPret} have shown that GRH implies (\dblast) with $C$ twice as large as predicted.
The Burgess bound seems even further from the truth: a folklore conjecture asserts that
$|S_\chi(t)| \ll_\e t^{1/2} q^\e$, from which it would immediately follow that
\begin{equation}
\tag{$\ddag$} \label{eq:VinogradovConj}
S_\xi(t) = o(t)
\qquad \forall t \gg_\epsilon p^{\epsilon} .
\end{equation}
This is known conditionally on GRH (or even on the weaker Generalized Lindel\"of Hypothesis).

Traditionally, P\'olya-Vinogradov and Burgess are considered to be somewhat independent of one another; the former is a global upper bound which applies only to long character sums, while the latter is a local bound which applies only to shortish sums.
Moreover, Burgess is often viewed as the superior result, both because in applications one frequently requires cancellation in short sums and because the proof is significantly more complex.
The goal of this note is to demonstrate that, in spite of appearances, the two bounds are intimately connected: even an apparently mild improvement of the `simpler' P\'olya-Vinogradov bound yields a significant improvement of the `deeper' Burgess bound.
More precisely, we'll show
\begin{BigThm}
Suppose the P\'olya-Vinogradov inequality \eqref{eq:PV} can be improved to $S_\chi(t) = o(\sqrt{q} \log q)$ for all even primitive quadratic $\chi \mod q$.
Then \eqref{eq:VinogradovConj}
holds for all odd primitive quadratic characters $\xi \mod p$.
\end{BigThm}

The proof of this theorem builds on joint work of the second author and Jonathan Bober \cite{BGnonsquare}, in which it is shown that any improvement of P\'{o}lya-Vinogradov yields an improvement of bounds on the least quadratic nonresidue. That theorem is quantitative, and shows that even improvements to the implicit constant in P\'olya-Vinogradov would break past the
$\frac{1}{4 \sqrt{e}}$
barrier in the least nonresidue problem.
In principle, it should be possible to make Theorem A quantitative as well.
To keep our exposition as brief and transparent as possible, we have elected not to do this, but it would be interesting to see what sorts of quantitative results one could obtain.

A crucial step in our argument is to show that if the mean value of a real multiplicative function up to some point is large, then the logarithmic mean up to that point must also be large.
Despite the extensive body of literature on multiplicative functions (or, perhaps, because of it!), we were unable to find a quantitative result of the form we need. To explain this further, we set some notation. Let
\[
\M_f(x) := \frac{1}{x} \sum_{n \leq x} f(n)
\qquad \text{and} \qquad
\L_f(x) := \frac{1}{\log x} \sum_{n \leq x} \frac{f(n)}{n}
\]
Theorem 2 of \cite{HalbRich} implies that for any multiplicative function $f : \Z \to [0,1]$,
\begin{equation}\label{eq:LogmeanLargerThanMean}
|\L_f(x)| \gg |\M_f(x)| .
\end{equation}
What if $f$ is allowed to take negative values?
It is a fun exercise to construct a family of completely multiplicative functions $f : \Z \to [-1,1]$ for which \eqref{eq:LogmeanLargerThanMean} fails.\footnote{For example, one can use Haselgrove's result \cite{Hasel} on the Liouville function to construct an infinite family of integers $N$ and a corresponding infinite family of completely multiplicative functions $f$ (depending on $N$) such that $\L_f(N) = 0 \neq \M_f(N)$.}
For our application, however, we will only need to compare the two means in the situation that both $x$ and $|\M_f(x)|$ are large, and in this case we will show that \eqref{eq:LogmeanLargerThanMean} does hold. More precisely:
\begin{MainLem}
Given $c > 0$, there exists $\delta = \delta(c) > 0$ and $x_0 = x_0(c) \geq 1$ such that
\[
|\M_f(x)| \geq c
\qquad \Longrightarrow \qquad
\L_f(x) \geq \delta
\]
for all completely multiplicative functions
$f : \Z \to [-1,1]$ and all $x > x_0$.
\end{MainLem}
\ni
It would be interesting to find a more explicit relationship between the two means which holds even when they are small. \\

\ni
\textit{Acknowledgements.}
The second author is grateful to a number of participants in the MSRI analytic number theory program for helpful discussions, particularly Adam Harper, Oleksiy Klurman, Dimitris Koukoulopoulos, and James Maynard. He would also like to thank Chantal David, K. Soundararajan, and Julia Wolf for generously providing office space. Last but not least, he would like to thank MSRI and the Simons Institute for the Theory of Computing for providing wonderful working environments for the mathematics community, where much of this work was carried out.

\section{Proofs} \label{sect:Proofs}

To streamline the proof of Theorem A, we isolate one of the key steps. Bober and the second author, building on previous work of the second author with Lamzouri \cite{GLevenchar,GLoddchar},
obtained omega results for character sums of a product of two characters.
In particular, the proof of Theorem 3 in \cite{BGnonsquare} (see equation (7) of that proof) implies:
\begin{lemma}\label{lem:BG}
Given odd primitive characters $\xi \mod k$ and $\psi \mod \ell$ such that $(k,\ell) = 1$. Consider the primitive character $\chi := \xi \psi$ of conductor
$q := k \ell$.
Then
\[
\frac{1}{\sqrt{q}} \max_{N \leq q} |S_{\chi}(N)|
\geq
\frac{\sqrt{\ell}}{\pi\varphi(\ell)}
\max_{t \leq q}
\Bigg|
\sum_{\substack{n \leq t \\ (n,\ell) = 1}} \frac{\xi(n)}{n}
\Bigg| + O(1) .
\]
\end{lemma}

\ni
We can now give a relatively short proof of Theorem A (assuming the validity of Lemma B, which we prove subsequently).

\begin{proof}[Proof of Theorem A]
Fix an $\e > 0$, and suppose the conjectured bound \eqref{eq:VinogradovConj} fails for some infinite collection of odd primitive real characters of prime conductor. More precisely, suppose there exists a positive constant $c$ and an infinite family
\[
\Xi := \{\text{primitive } \; \xi \mod p :
    \xi(-1) = -1 , \xi^2 = \chi_0\}
\]
such that for each $\xi \in \Xi$ we have
\begin{equation} \label{eq:MeanValOrigForm}
|S_\xi(t_p)| \geq c t_p
\end{equation}
for some $t_p > p^\e$.
We will construct an infinite family of even primitive real characters $\chi\mod q$ satisfying
\begin{equation} \label{eq:GoalOfPrfA}
\max_{N} |S_\chi(N)| \gg_{c,\e} \sqrt{q} \log q ,
\end{equation}
thus contradicting the hypothesis of Theorem A.

Reformulating \eqref{eq:MeanValOrigForm}
to read
\(
|\M_\xi(t_p)| \geq c ,
\)
we are led to apply Lemma B.
The lemma produces two positive constants $\d$ and $x_0$, both depending only on $c$.
We may assume (after possibly removing a finite number of characters from $\Xi$) that $p^\e > x_0$ for all $\xi \mod p \in \Xi$.
Thus $t_p > x_0$, whence Lemma B implies
\begin{equation}\label{eq:LowerBdOnLogAvg}
\L_\xi(t_p) \geq \d .
\end{equation}
Let $\ell$ be the smallest prime larger than $\frac{2}{\d}$ which satisfies $\ell \equiv 3 \mod 4$, and set
\[
\psi := \leg{\cdot}{\ell} .
\]
To each character $\xi \mod p \in \Xi$ we associate a character $\chi \mod q$ defined by
\[
\chi := \xi \psi .
\]
Note that each such $\chi$ is an even primitive real character with conductor $q = p \ell$, whence Lemma \ref{lem:BG} implies
\begin{equation}\label{eq:DesiredIneqForA}
\frac{1}{\sqrt{q}} \max_{N \leq q} |S_{\chi}(N)|
\gg_{c,\e}
\Bigg|
\sum_{\substack{n \leq t_p \\ \ell \dnd n}} \frac{\xi(n)}{n}
\Bigg| + O(1) .
\end{equation}
We now show that the right hand side is $\gg_{c,\e} \log q$. It may be helpful to recall the dependencies among our parameters: $c$ and $\e$ are fixed, $\d$ depends only on $c$, and $\ell$ depends only on $\d$ (and hence only on $c$).
We have
\[
\begin{split}
\sum_{\substack{n \leq t_p \\ \ell \dnd n}} \frac{\xi(n)}{n}
& =
\sum_{n \leq t_p} \frac{\xi(n)}{n} -
\frac{\xi(\ell)}{\ell}
    \sum_{m \leq t_p / \ell} \frac{\xi(m)}{m} \\
& \geq
\L_\xi(t_p) \log t_p -
\frac{1}{\ell}
    \Big(\log (t_p / \ell) + \gamma + O(\ell/t_p)\Big) \\
& \geq
(\d - \frac{1}{\ell}) \log t_p + \frac{\log \ell - \gamma}{\ell} + O(1/t_p) \\
& \geq
\frac{\d \e}{2} \log p + O(1/p^\e) \\
& =
\frac{\d \e}{2} \log q - \frac{\d \e}{2} \log \ell + O(1/p^\e)
\end{split}
\]
where we have used \eqref{eq:LowerBdOnLogAvg} combined with $\ell > \frac{2}{\d} > 2$ and $t_p > p^\e$.
Substituting this into \eqref{eq:DesiredIneqForA} yields \eqref{eq:GoalOfPrfA} as desired.
\end{proof}

\ni
Having proved Theorem A under the assumption of Lemma B, it therefore suffices to handle the lemma.

\begin{proof}[Proof of Lemma B]
Fix $c > 0$, and suppose
\[
|\M_f(x)| \geq c
\]
for some completely multiplicative function $f : \Z \to [-1,1]$ and some large fixed $x$. Our goal is to show that
\[
\L_f(x) \gg_c 1 .
\]
The key subtlety here is that the $x$ appearing in the previous two displays is the same.

Our first step is classical -- we approximate the logarithmic mean of $f$ by the mean of $\one * f$, where $\one$ denotes the constant function taking the value 1 and $*$ denotes the Dirichlet convolution:
\begin{equation}\label{eq:DirichletConvolutionApprox}
\sum_{n \leq x} \frac{f(n)}{n} =
\frac{1}{x} \sum_{n \leq x} (\one * f)(n) + O(1) .
\end{equation}
This reduces our problem to bounding from below the right hand side of the above. There are many ways to do this, but we will take a shortcut and simply quote from Granville-Soundararajan's study \cite{GS} of the minimum of the logarithmic mean of a multiplicative function. Equation (3.5) of that paper implies
\begin{equation}\label{eq:GSLowerBd}
\frac{1}{x} \sum_{n \leq x} (\one * f)(n) \gg
e^{-u e^{u/2}} \log x + O(1) ,
\end{equation}
where
\begin{equation}\label{eq:DefnOfu}
u := \sum_{p \leq x} \frac{1 - f(p)}{p} .
\end{equation}
Bounding this from below amounts to bounding $u$ from above. We will deduce such a bound from the Hall-Tenenbaum theorem on mean values of multiplicative functions, although many other theorems of similar flavor would also suffice. A special case of the main theorem of \cite{HallTen} asserts the existence of a constant $\kappa \approx 0.32$ such that
\[
|\M_f(x)| \ll e^{-\kappa u} ,
\]
where $u$ is defined by \eqref{eq:DefnOfu}.
On the other hand, we are assuming
$|\M_f(x)| \geq c$, whence
\[
u \ll_c 1 .
\]
Plugging this into \eqref{eq:GSLowerBd} and using this in the estimate \eqref{eq:DirichletConvolutionApprox} produces the bound
\[
\sum_{n \leq x} \frac{f(n)}{n} + O(1)
    \gg_c \log x + O(1)
\]
where the two $O(1)$ terms are bounded by some constant which is independent of $c$. Taking $x$ large enough (where `large' depends only on $c$) we can make the contribution of the $O(1)$ terms negligible. This concludes the proof.
\end{proof}


\end{document}